\documentclass[onecolumn]{cdcarticle}
\usepackage{amsfonts}
\usepackage{amssymb}
\usepackage{amsmath}
\usepackage{amsthm}
\usepackage{graphicx}

\usepackage{stdmath}
\usepackage{overpic}

\graphicspath{{graphics_included/}}

\def\lcplus{\textbf{ edited}}

\def\lcversion{2006.03.21.01}

\def\onevec{\mathbf{1}}

\begin{document}

\title{Degree Bounds for Polynomial Verification of\\
  the Matrix Cube Problem}

\author{%
  Been-Der Chen\footnotesymbol{1,3}
  \and
  Sanjay Lall\footnotesymbol{2,3}
  }

\note{Submitted to Mathematical Programming}

\maketitle

\makefootnote{1}{Email: bdchen@stanford.edu}

\makefootnote{2}{Email: lall@stanford.edu}

\makefootnote{3}{%
  Department of Aeronautics and Astronautics,
  Stanford University, Stanford CA 94305-4035, U.S.A.}

\makefootnote{}{
  Both authors were partially supported by the Stanford MICA
  \emph{Coordination and  Control for Networks of Interacting
    Automata}, DARPA/SPAWAR award number N66001-01-C-8080
  }


\begin{abstract}
  In this paper we consider the problem of how to computationally test
  whether a matrix inequality is positive semidefinite on a
  semialgebraic set.  We propose a family of sufficient conditions
  using the theory of matrix Positivstellensatz refutations.  When the
  semialgebraic set is a hypercube, we give bounds on the degree of
  the required certificate polynomials.
\end{abstract}


\nocite{*}

\section{Introduction}

In this paper we consider the following problem.

\begin{problem}
  \label{prob:main}
  Suppose $H_0,\dotsc,H_m\in\R^{n\times n}$ are symmetric matrices,
  and $\Delta\subset\R^m$. Define the affine map
  $G:\R^m\rightarrow \R^{n\times n}$ by 
  \[
     G(\delta) = H_0 + \sum_{i=1}^m \delta_i H_i
  \]
  for all $\delta \in \R^m$. We would like to know if 
  \[
  G(\delta) \geq 0 \text{ for all }\delta \in\Delta.
  \]
\end{problem}
This problem is a \emph{robust semidefinite program}, and it has many
important applications in control and optimization. One motivating
application is testing \emph{quadratic stability}, as follows.
Consider the parameterized family of linear time-invariant systems
\begin{align*}
  \dot x &= \bbl(A_0 + \sum_{i=1}^m \delta_i A_i \bbr) x
\end{align*}
Here $\delta\in\R^m$ is a vector of \emph{uncertain parameters}.  We
would like to check whether the above system is stable for all
$\delta\in\Delta$. This problem has been addressed
in~\cite{Popov_1960, YOUNG_2001, SAF_1986}. A well-known approach is
to seek a quadratic Lyapunov function which proves stability for all
parameters within the uncertainty set $\Delta$. That is, we
would like to find a positive definite matrix $P\in\R^{n\times n}$ such that
\begin{align*}
  \bbl(A_0+\sum_{i=1}^m \delta_i A_i\bbr)^T P + 
  P \bbl(A_0+\sum_{i=1}^m \delta_i A_i\bbr) < 0
\end{align*}
for all $\delta\in\Delta$. Testing whether $P$ satisfies this
inequality is equivalent to Problem~\ref{prob:main}, via the
identification $H_i=-A_i^T P- P A_i$ for $i=0,\dotsc,m$.

More generally, we can convert a more general class of robust
optimization problems to the form of Problem~\ref{prob:main}. We would
like to solve
\begin{align*}
  \min \qquad & c^T x \\
  \text{subject to } \quad & \mathcal{A}_0 + 
     \sum_{i=1}^m x_i \mathcal{A}_i \geq 0  
     \qquad \text{for all } 
     (\mathcal{A}_0,\dotsc,\mathcal{A}_n)\in\mathcal{U}
\end{align*}
where the set of matrix tuples $\mathcal{U}$ is given by 
\begin{equation*}
  \mathcal{U} = \bbl\{
       ( A_1,\dotsc,A_n )= \hskip 32mm \\
         (A_0^0,\dotsc,A_n^0) + \sum_{k=1}^m \delta_i (A_0^k,\dotsc,A_n^k)
    \,\verttt \,\, \delta\in\Delta
  \bbr\}
\end{equation*}
To find the optimal solution to this robust semidefinite problem, we
need to be able to efficiently verify that a given $x$ satisfies the
constraints, and this is equivalent to Problem~\ref{prob:main}. In
fact, even this verification problem is hard for most uncertainty
sets.

Problem~\ref{prob:main} has been addressed in the literature.  When
the uncertainty set is an ellipsoid, the problem may be converted to a
binary optimization problem~\cite{Bental_Nemirovski_1998}.  When the
uncertainty set is a hypercube, Problem~\ref{prob:main} is called the
\emph{matrix cube problem}, and it was shown to be NP hard
in~\cite{Nemirovski_93}. In the case when the uncertainty set is a
bounded polytope, it is sufficient to check the matrix inequality at
the vertices. Notice that in the case when $\Delta$ is a cube there
are $2^m$ vertices, and so this approach scales very poorly as $m$
grows. Similar results for the quadratic stability problem are shown
in~\cite{HG_1976},~\cite{Boyd_Young_1989}.

To reduce computational complexity , several sufficient conditions
have been proposed, such as the use of the $\mathcal{S}$-procedure to
construct a set of scalar certificates~\cite{BGF_1994}. Ben-Tal and
Nemirovski also proposed a stronger condition which does not exhibit
the above poor scaling~\cite{BN_2002}. In this paper we will
generalize this condition, so we state it here. Here $\S^n$ denotes
the set of real $n\times n$ symmetric matrices.

\begin{thm}
\label{thm:thm_bental}
Suppose $\Delta$ is the cube
\[
\Delta=\bl\{\,
\delta\in\R^m \ \vert \ 
\abs{\delta_i} \leq 1 \text{ for all }i
\,\br\}
\]
Define the set $\mathcal{X}_T\subset\S^n \times \dots \times \S^n$,
where $(X_1,\dotsc,X_m)\in\mathcal{X}_T$ if and only if
\begin{align*}
  H_0 - \sum_{i=1}^m X_i &\geq 0, \\
  X_i + H_i &\geq 0, \qquad \text{for all } i=1,\dotsc,m\\
  X_i - H_i &\geq 0, \qquad \text{for all } i=1,\dotsc,m.
\end{align*}
Then $G(\delta)$ is positive semidefinite for all $\delta\in\Delta$
if $\mathcal{X}_T$ is not empty.
\end{thm}

This condition may be tested via semidefinite programming.  The
paper~\cite{BN_2002} also shows that if the above SDP condition is
infeasible then there exists a $\delta$ within a larger cube such that
$G(\delta)$ is not positive semidefinite. This gives an estimate of
the conservatism of this test.

The matrix cube problem (and also therefore Problem~\ref{prob:main})
is closely related to binary quadratic programming. Here, we would
like to find
\begin{align*}
  \min \qquad & x^T A x \\
  \text{subject to } \quad & x\in\{-1, 1\}^n
\end{align*}
Without loss of generality we may assume $A$ is positive definite,
and it is then straightforward to see that the problem
is equivalent to the following matrix cube problem.
\begin{align*}
  \max \qquad & t \\
  \text{subject to } \quad &
  \bmat{ t & x^T \\ x & A^{-1} } \geq 0 
            \quad \text{for all } x\in\{-1,1\}^n
\end{align*}

In general such quadratic programs are hard, and much research has
been done to address this, for example using the Lagrangian relaxation
to compute a lower bound on the optimal value~\cite{Shor_1987}, or
using semidefinite programming via a lifting
approach~\cite{Lovasz_Schrijver_1991}. The gap between the relaxed
problem and the actual problem may be reduced by introducing
additional variables and redundant
constraints~\cite{Anstreicher_Wolkowicz_2000}.  Lasserre used an
approach based on moments and showed that one needs at most $2^m-1$
additional variables~\cite{Lasserre_2001, Lasserre_01_2001} for an
exact solution.  This approach is also related to the dual of the
refutation approach adopted in this paper.

In this paper, instead of searching a set of scalar certificates using
the $\mathcal{S}$-procedure, we will construct a sufficient condition
via a search for a polynomial certificates.  If such a certificate
exists, then there is no $\delta\in\Delta$ such that the affine
function $G$ is not positive semidefinite.  Our approach is applicable
to general semi-algebraic uncertainty sets, including ellipsoids and
hypercubes. In this formulation, we construct a family of refutation
sets which have a hierarchical structure. If the current refutation
set does not yield a feasible certificate, we may seek for higher
degree certificates.  Similar approaches have been used to analyze and
synthesize output feedback controllers for LPV systems~\cite{WP_2004}.

For some uncertainty sets we will also show that if there is no
$\delta\in\Delta$ for which $G$ is not positive semidefinite, then
there will exists a certificate of specific degree.  When the
uncertainty set is a hypercube, we show that the highest degree needed
is at most $2m$. We also study the case when the certificates are
restricted to be quadratic and we show the resulting condition is
tighter than the best existing result of Theorem~\ref{thm:thm_bental}.
In addition, we give several cases when our conditions using quadratic
certificates are necessary and sufficient. Finally, we give some
numerical examples to compare our results with others.


\section{Preliminaries}

We use the following standard notation. The matrix $I_n$ is the $n
\times n$ identity. For $X\in\S^n$, the notation $X>0$ means that $X$
is positive definite. The vector $e_i\in\R^n$ has the $i$th entry
equal to 1 and all other entries equal to zero. The vector
$\onevec\in\R^n$ has all entries equal to 1.

The set $\R[x_1,\dots,x_n]$ is the ring of polynomials in $n$
variables with real coefficients. We often abbreviate
$\R[x_1,\dots,x_n]$ to simply $\R[x]$. Every polynomial
$f\in\R[x]$ can be written as
\[
f= \sum_{\alpha\in W} c_{\alpha} x^{\alpha}
\]
where $W \subset \N^n$, and the notation $x^\alpha$ is defined by
\[
x^\alpha = x_1^{\alpha_1}  x_2^{\alpha_2} \dots   x_n ^{\alpha_n}
\]
A polynomial $g\in\R[x]$ is called a \emph{sum of squares (SOS)} if it
can be expressed as
\[
  g(x)= \sum_{i=1}^n f_i(x)^2
\]
for some polynomials $f_i\in\R[x]$. We use $\Sigma[x]$ to represent
the set of sum-of-squares polynomials in $\R[x]$, and abbreviate it to
$\Sigma$ when the dimension is clear from the context. We also extend
this definition matrix polynomials as follows.. Let $\R[x]^{m\times
  n}$ denote the set of $m\times n$ polynomial matrices and $\S[x]^n$
denote the set of $n\times n$ symmetric polynomial matrices. We define
the notion of \emph{sum-of-squares} for matrix polynomials as follows
\begin{defn}
  A matrix polynomial $S\in\S[x]^{m}$ is called a sum-of-squares if
  there exist a matrix polynomial $T\in\R[x]^{m\times q}$ such that
  \[
     S(x) = T(x) T(x)^T.
  \]
\end{defn}
This is a generalization of SOS representation used for scalars.  We
will use $\Sigma[x]^n$ to represent the set of $n\times n$ SOS
polynomial matrices. We also define two specific sets $\mathcal{Q}_1$,
$\mathcal{Q}_2$ which will be useful in later sections.
\begin{defn}
  Let $W_1=\{\alpha\in\N^m\,|\,\alpha_i\leq 2 \text{ for all
  }i=1,\dotsc,m\}$ and $W_2=\{\alpha\in\N^m\,|\,\sum_{i=1}^m
  \alpha_i\leq 2 \}$. The sets $\mathcal{Q}_1$, $\mathcal{Q}_2$ are
  defined as
  \begin{align*}
    \mathcal{Q}_1 &= \bbl\{ \sum_{\alpha\in W_1} C_{\alpha}
    \delta^{\alpha} \,|\, C_{\alpha}\in\S^n \text{ for all } \alpha\in
    W_1 \bbr\}
    \\
    \mathcal{Q}_2 &= \bbl\{ \sum_{\alpha\in W_2} C_{\alpha}
    \delta^{\alpha} \,|\, C_{\alpha}\in\S^n \text{ for all } \alpha\in
    W_2 \bbr\}
  \end{align*}
\end{defn}
Note that polynomials in $\mathcal{Q}_1$ have degree less than or
equal to $2m$ and polynomials in $\mathcal{Q}_2$ have degree less than
or equal to $2$. 

When $F\in\Sigma[x]^n$, it is clear that $F$ is positive semidefinite
for all $x\in\R^n$. One may address positive semidefiniteness of a
matrix polynomial over a restricted domain using the following lemma,
which gives a sufficient condition.

\begin{lem}
  \label{thm:putinar} Suppose $f_1,\dotsc,f_n\in\R[x]$ and
  $Q\in\R[x]^m$ is a symmetric matrix polynomial.  Define the set
  \[
  \mathcal{D}= \bbl\{\, x\in\R^n \ \vert \ f_i(x)\geq 0 \text{ for all }
  i=1,\dotsc,m \,\bbr\}
  \]
  Then, $Q(x) \geq 0$ for all $x\in\mathcal{D}$ if there exist SOS
  polynomial matrices $S_0,S_1,\dotsc,S_n\in\Sigma[x]^m$ such that
  \begin{align*}
    Q(x) &= S_0(x) + \sum_{i=1}^n S_i(x) f_i(x)
  \end{align*}
\end{lem}
It is also known that if $\mathcal{D}$ is compact, and with additional
technical restrictions, then the above condition is also
necessary~\cite{Scherer_Hol_2004}; this is an extension of a
well-known result by Putinar~\cite{Putinar_1993}.

\section{Positivstellensatz refutations}

In this section, we will study Problem~\ref{prob:main} when the set
$\Delta$ is semialgebraic, that is
\begin{align*}
  \Delta &= \bbl\{ \delta\in\R^m \,|\,  f_i(\delta)\geq 0,
                \text{ for } i=1,\dotsc,m \bbr\}
\end{align*}
where $f_1,\dotsc,f_m\in\R[\delta]$. It is clear that a cube and an
ellipsoid can be expressed as semi-algebraic sets. The following
condition provides a simple condition under which $G(\delta)$ is
positive semidefinite for all $\delta\in\Delta$.

\begin{thm}
  \label{thm:thm_positive}
  The matrix polynomial $G(\delta)$ is positive semidefinite for all
  $\delta\in\Delta$ if there exist SOS polynomial matrices
  $S_0,S_1,\dotsc,S_m$ satisfying
  \begin{equation}
    \label{thm:refutation}
     G(\delta) =  S_0 + \sum_{i=1}^m S_i f_i(\delta)
  \end{equation}
\end{thm}

This is a simple consequence of Lemma~\ref{thm:putinar} and we may
view it as provided a certificate refuting the existence of
$\delta\in\Delta$ such that $G(\delta)$ is not within the positive
semidefinite cone. The certificate is the sequence of polynomials
$S_0,\dots,S_m$. As discussed in the previous section, this condition
is also necessary if additional technical conditions on $\Delta$ are
satisfied~\cite{Scherer_Hol_2004} (both polytopes and ellipsoids
satisfy these conditions.)

One thing we have not yet specified is the degree of the certificates
required.  Although we may pursue high degree certificates, the
computational complexity of finding $S_0,\dots,S_m$ grows rapidly as
we search over sets containing high-degree polynomials.  In many
applications of this refutation approach, a bound on the degree of the
required certificates is not known.  However, in some cases, we can
show a degree bound. First, we consider the case when the uncertainty
set $\Delta$ is a simplex.

\begin{thm}
  Suppose $\Delta$ is the simplex 
  \begin{align*}
    \Delta=\{\delta\in\R^m\,|\,\delta\geq 0,\onevec^T \delta\leq 1\}
  \end{align*}
  Then $G(\delta)\geq 0$ for all $\delta\in\Delta$ if and only if
  there exist positive semidefinite matrices
  $S_0,\dotsc,S_{m+1}\in\S^n$ such that
  \begin{align}
    \label{eqn:refutation_convex_hull}
    G(\delta) &= S_0 + 
       \sum_{i=1}^m \delta_i S_i + (1-\onevec^T \delta) S_{m+1}
  \end{align}  
\end{thm}

\begin{proof}
  Sufficiency is straightforward as in Theorem~\ref{thm:thm_positive}
  and we will now prove necessity. Assume that $G(\delta)\geq 0$ for
  all $\delta\in\Delta$, then
  \begin{align*}
    H_0 &\geq 0 \\ H_0 + H_i &\geq 0 \qquad \text{for } i=1,\dotsc,m
  \end{align*}
  It can be verified that $S_0=0$, $S_i=H_0+H_i$ for $i=1,\dotsc,m$
  and $S_{m+1}=H_0$ satisfy~\eqref{eqn:refutation_convex_hull} and
  they are positive semidefinite.
\end{proof}

The above theorem precisely specifies the certificates in the very
special case that the the uncertainty set is the simplex. If the
uncertainty set is expressed as a convex hull of the vertices, one can
also convert the verification problem to the above form and apply the
same refutation.

Although the theorem gives necessary and sufficient conditions, the
number of required certificates is the same as the number of vertices.
For the hypercube which has $2^m$ vertices, this soon becomes
intractable as $m$ grows. On the other hand, the hypercube has another
representation as a semialgebraic set
\[
\Delta
=\{\delta\in\R^m\,|\,\delta_i^2\leq 1,\text{for }i=1,\dotsc,m\}
\]
and in the following result we use this representation.

\begin{thm}
  Define the set $\mathcal{X}_1\subset\Sigma[\delta]^n\times
  \S[\delta]^{n} \times \dots \times\S[\delta]^{n}$ such that
  $(S_0,S_1,\dotsc,S_m)\in\mathcal{X}_1$ if and only if 
  \begin{align}
    \label{eqn:cube_refutation}
    G(\delta) &= S_0 + \sum_{i=1}^m (1-\delta_i^2) S_i.
  \end{align}
  Then, $G(\delta)\geq 0$ for all $\delta$ within the unit cube if and
  only if $\mathcal{X}_1$ is not empty.
\end{thm}

\begin{proof}
  If there exist $(S_0,\dotsc,S_m)$
  satisfying~\eqref{eqn:cube_refutation}, then it is clear that the
  right hand side of~\eqref{eqn:cube_refutation} is positive
  semidefinite for all $\delta$ at the vertices of the cube, and this
  implies that $G(\delta) \geq 0$ for all $\delta\in\Delta$.

  Now, we show the converse. The set of vertices may also be
  represented as the set of $\delta$ satisfying $1-\delta_i^2\geq 0$
  and $\delta_i^2-1\geq 0$ for $i=1,\dotsc,m$. By choosing
  $r=\sqrt{m}$, we can verify that
  \begin{align*}
    r^2 -\norm{\delta}^2 &= 2\sum_{i=1}^m (1-\delta_i^2) + 
        \sum_{i=1}^m (\delta_i^2-1).
  \end{align*}
  Since $G(\delta) \geq 0$ for all $\delta$ at the vertices, there
  exist SOS matrix polynomials $P_0, P_1,\dotsc, P_m$,
  $P_{m+1},\dotsc,P_{2m}$ such that
  \begin{align*}
    G(\delta) &= P_0 + \sum_{i=1}^m (1-\delta_i^2) P_i
            + \sum_{i=1}^m (\delta_i^2-1) P_{m+i}
  \end{align*}
  from Theorem 2 in~\cite{Scherer_Hol_2004}. Letting $S_0=P_0$ and
  $S_i=P_i-P_{m+i}$ for $i=1,\dotsc,m$ completes the proof.

\end{proof}
The above theorem shows that $S_1,\dotsc,S_m$ only need to be
symmetric matrix polynomials instead of SOS matrix polynomials. Now,
we show the refutation is degree bounded and the highest degree of
certificates required is $2m$.

\begin{thm}
  \label{thm:matrix_cube}
  Define the set $\mathcal{X}_2 \subset \mathcal{X}_1$ such that
  $(S_0,\dotsc,S_m)\in\mathcal{X}_2$ if and only if
  $(S_0,\dotsc,S_m)\in\mathcal{X}_1$ and $S_1,\dotsc,S_m \in
  \mathcal{Q}_1$. Then, $G(\delta)\geq 0$ for all $\delta$ within the
  unit cube if and only if $\mathcal{X}_2$ is non-empty.
\end{thm}

Before proving the theorem, we show two lemmas.

\begin{lem}
\label{lem:LMI_stack} 
Suppose there are two symmetric matrices $A$,$B\in\S^n$ satisfying $-A
\leq B $ and $B \leq A$. Then
  \begin{align*}
    \bmat{A & B \\ B & A} \geq 0.
  \end{align*}
\end{lem}

\begin{proof}
  Suppose $-A \leq B$ and $B \leq A$. Then
  \begin{align*}
    \bmat{A & B \\ B& A} &= \frac{1}{2} \bmat{I & I \\ -I & I}
                     \bmat{A-B \\ & A+B} \bmat{I & I \\ -I & I}^T \geq 0.
  \end{align*}
\end{proof}

\begin{lem}
  \label{lem:matrix_cube_LMI}
  Suppose $H_0,\dotsc,H_m\in\S^n$.  Define the block diagonal matrices
  \[
  G_k = \bmat{H_k \\ & \ddots \\ & & H_k}\in\S^{2^{k-1}n},
              \quad \text{for }k=1,\dots,m.
  \]
  and recursively define the sequence $N_0,N_1,\dots,N_m$ by
  \[
  N_k = \bmat{N_{k-1} & G_k \\ G_k & N_{k-1}}
  \qquad
  \qquad
  N_0=H_0
  \]
  Then $N_m\geq0$ if
  \begin{align}
  \label{eqn:eq_matrix_cube}
    H_0 + \sum_{i=1}^m \delta_i H_i \geq0, \qquad
                        \text{for all } \delta\in\{-1,1\}^m.
  \end{align}
\end{lem}

\begin{proof}
  We prove this by induction on $m$. First notice that the result is
  trivially true when $m=0$.  To show the case when $m=k+1$, we have
  \begin{align*}
    H_0 + \sum_{i=1}^{k+1} \delta_i H_i \geq0
    \qquad \text{for all } \delta_1,\dots,\delta_{k+1}\in\{-1,1\}
  \end{align*}
  if and only if
  \begin{align*}
    (H_0+H_{k+1}) + \sum_{i=1}^k \delta_i H_i &\geq 0, &
    (H_0-H_{k+1}) + \sum_{i=1}^k \delta_i H_i &\geq 0
  \end{align*}
  holds whenever $\delta_1,\dots,\delta_k \in \{-1,1\}$.  Now assume
  the lemma holds for $m=k$, and this then implies
  \begin{align*}
    \bmat{N_{k-1} & G_k\\ G_k& N_{k-1}} + G_{k+1} &\geq 0, &
    \bmat{N_{k-1} & G_k\\ G_k& N_{k-1}} - G_{k+1} &\geq 0.
  \end{align*}
  Applying Lemma~\ref{lem:LMI_stack} proves the theorem.
\end{proof}

Now, we prove Theorem~\ref{thm:matrix_cube}.

\medskip {\noindent\bf Proof of Theorem~\ref{thm:matrix_cube}.}
Sufficiency is implied by Theorem~\ref{thm:thm_positive} and we will
now prove necessity. Suppose $G(\delta)\geq 0$ for all $\delta$ within
the unit cube. It is shown in Lemma~\ref{lem:matrix_cube_LMI} that
$N_m \geq 0$. Let $z_0=1$ and define a sequence of vectors of
monomials $z_0,z_1,\dots,z_m$ recursively as follows
  \begin{align*}
    z_{i} &= \bmat{z_{i-1} \\ \delta_{i} z_{i-1}}
                        \qquad \text{for } i=1,\dotsc,m.
  \end{align*}
  Choose $S_0$ as
  \begin{align*}
    S_0 = 2^{-m} (z_m\otimes I)^T N_m (z_m\otimes I).
  \end{align*}
  It is clear that $S_0\in\Sigma[x,\delta]\cap\mathcal{Q}_1$.

  As for $S_1,\dotsc,S_m$, we let $\mathcal{D}=\{1,\dotsc,m\}$ and
  define the sets $\mathcal{E}_{k,l},\mathcal{F}_{j,k,l}$ as follows,
  \begin{align*}
    \mathcal{E}_{k,l} &= \bbl\{\, A\subset\mathcal{D}
      \ \vert \
      |A|=l \text{ and } k\not\in A \,\bbr\}, \\
    \mathcal{F}_{j,k,l} &= \bbl\{\, A\subset\mathcal{D}
      \ \vert \
      |A|=l, \text{ and } j,k\not\in A \,\bbr\}.
  \end{align*}
  We also define the polynomials
  \begin{align*}
    p_{k,l} &= \sum_{A\in \mathcal{E}_{k,l}} \prod_{j\in A} \delta_j^2, &
    q_{j,k,l} &= \sum_{A\in \mathcal{F}_{j,k,l}} \prod_{p\in A} \delta_p^2
  \end{align*}
  Let $S_1,\dotsc,S_m$ as follows
  \begin{align*}
      S_k = H_0 \left(c_1+\sum_{i=1}^{m-1} c_{i+1} p_{k,i}\right)
       +\sum_{\substack{i=1\\i\ne k}}^m \delta_i H_i
     \left(d_1+\sum_{j=1}^{m-2} d_{j+1} q_{i,k,j}\right)
      \\
     \quad \text{for } k=1,\dotsc,m
  \end{align*}
  where $c\in\R^m$ and $d\in\R^{m-1}$ satisfy
  \begin{align*}
    M(m) c &= e_1-2^{-m}\onevec, &
    M(m-1) d &= e_1-2^{-m+1}\onevec
  \end{align*}
  and $M:\R \mapsto \R^{m \times m}$ are
  \begin{align*}
     M(m)&= \bmat{m & 0 & \cdots & 0 \\ -1 & m-1 & \ddots & \vdots \\
        & \ddots & \ddots & 0 \\ 0 & \cdots & -m+1 & 1}.
  \end{align*}
  The highest degree of $S_1,\dotsc,S_m$ in each $\delta_i$ is at most
  2 which implies that $S_1\dotsc,S_m\in\mathcal{Q}_1$. Expanding
  $S_0+\sum_{i=1}^m (1-\delta_i^2)S_i$ shows that $S_0,\dotsc,S_m$
  satisfy~\eqref{eqn:cube_refutation}.
\hfill\qed

The reason that degree bounded $S_0, S_1,\dotsc,S_m$ exist is because
of the persymmetric structure of $N_m$. We now have a family of
refutations for the matrix cube problem. If $G(\delta)\geq 0$ for all
$\delta\in\Delta$, then we have shown that there will exist a
certificate of degree at most $2m$.

This condition may also be directly expressed as a semidefinite
program.  Although the degree bound on the certificates grows linearly
with $m$, the number of monomials required to express the certificates
(i.e., the dimension of $\mathcal{Q}_1$) grows exponentially in $m$.
Of course, since we have exactly solved the problem this is expected;
the original problem is NP-hard.

When the cube is high-dimensional, the computational complexity of
searching $\mathcal{X}_2$ is also high.  To reduce computational
effort, we may limit the search to low degree certificates. We now
consider this case, limiting the search to the set $\mathcal{X}_3$ as
follows.

\begin{defn}
  Define the set $\mathcal{X}_3 \subset \mathcal{X}_2$ such that
  $(S_0, \dotsc, S_m) \in \mathcal{X}_3$ if and only if
  $(S_0,\dotsc,S_m)\in\mathcal{X}_2$ and $S_0\in\mathcal{Q}_2$,
  $S_1,\dotsc,S_m\in\S^n$.
\end{defn}

Testing if $\mathcal{X}_3$ is nonempty is equivalent to the following
semidefinite program.
\begin{equation}
  \label{eqn:quadratic_refutation}
  \begin{aligned}
    \text{find } \qquad & X_1,\dotsc,X_m\in\S^n, L\in\S^{n(m+1)}
        \\[1mm]
    \text{such that } & 
    L = \bmat{L_{00} & \cdots & L_{0m} \\
                \vdots & \ddots & \vdots \\
                L_{0m}^T & \cdots & L_{mm}} \geq 0 \\
       & 0 = \sum_{i=0}^m L_{ii} - H_0 \\
       & 0 = L_{ii} - X_i \hskip 20mm \text{for } i=1,\dotsc,m \\
       & 0 = L_{0i} + L_{0i}^T - H_i \hskip 11mm \text{for } i=1,\dotsc,m \\
       & 0 = L_{ij} + L_{ij}^T \qquad \text{for } i,j=1,\dotsc,m,i\ne j
  \end{aligned}
\end{equation}
If the above semidefinite program is feasible, we may choose $S_i=X_i$
for $i=1,\dotsc,m$ respectively. This implies that $S_1,\dotsc,S_m$
are positive semidefinite.

The gap between verifying the matrix cube problem and checking the
non-emptiness of $\mathcal{X}_3$ can be interpreted as the degree
difference between certificates in $\mathcal{X}_2$ and
$\mathcal{X}_3$. The degree of certificates in $\mathcal{X}_3$ is at
most $2$, instead of growing linearly with respect to $m$.  Although
this means that the condition is conservative, we now show that it is
still tighter than the previously well-known condition in
Theorem~\ref{thm:thm_bental}.

\begin{thm}
  If $\mathcal{X}_T$ is not empty, then $\mathcal{X}_3$ is not empty.
\end{thm}

\begin{proof}
  Suppose $(X_1,\dotsc,X_m)\in\mathcal{X}_T$. It is clear that $X_i$
  is positive definite for $i=1,\dotsc,m$. From
  Lemma~\ref{lem:LMI_stack}, $X_1,\dotsc,X_m$ also satisfy
  \begin{align*}
    \bmat{X_i & H_i \\ H_i & X_i}\geq 0, \qquad \text{for } i=1,\dotsc,m
  \end{align*}
  which by the Schur complement implies
  \begin{align*}
    X_i \geq H_i X_i^{-1} H_i, \qquad \text{for } i=1,\dotsc,m.
  \end{align*}
  Thus,
  \begin{align*}
    H_0 - \frac{1}{2}\sum_{i=1}^m X_i -
              \frac{1}{2} \sum_{i=1}^m H_i X^{-1} H_i &\geq
    H_0 - \sum_{i=1}^m  X_i \\ &\geq 0.
  \end{align*}
  Applying the Schur complement again gives
  \begin{align}
    \label{eqn:LMI_bental}
    J &= \bmat{H_0-\frac{1}{2} \sum_{i=1}^m X_i &
            \frac{1}{2} H_1 & \cdots & \frac{1}{2} H_m \\[1mm]
          \frac{1}{2} H_1 & \frac{1}{2} X_1 & & 0 \\[1mm]
          \vdots & & \ddots \\[1mm]
          \frac{1}{2} H_m & 0 & & \frac{1}{2} X_m} \geq  0
  \end{align}
  Let
  \begin{align*}
    S_0 &= \bmat{I \\ \delta\otimes I}^T J
                     \bmat{I \\ \delta\otimes I}, &
    S_i &= \frac{1}{2} X_i, \quad \text{for } i=1,\dotsc,m.
  \end{align*}
  It is clear that $S_1,\dotsc,S_m$ are positive semidefinite. Also
  $S_0\in\Sigma\cap\mathcal{Q}_2$.  Expanding $S_0$ shows
  $(S_0,\dotsc,S_m)\in\mathcal{X}_3$.
\end{proof}

The above theorem shows that every certificate in $\mathcal{X}_T$ has
a corresponding instance in $\mathcal{X}_3$. We show that our
condition is strictly tighter than the previous condition via a
counterexample in Section~\ref{sec:exampleone}.  Comparing the two
semidefinite programs, we recognize that~\eqref{eqn:LMI_bental}
imposes constraints on the off-diagonal entries such that $L_{ij}=0$
for $i,j=1,\dotsc,m,i\ne j$.  The entries are relaxed to be skew
symmetric in~\eqref{eqn:quadratic_refutation} and this condition is
still sufficient.  This skew symmetric structure arises naturally in
the Positivstellensatz refutation.

\begin{figure}[htb]
\centerline{
  \begin{overpic}[width=0.55\textwidth]{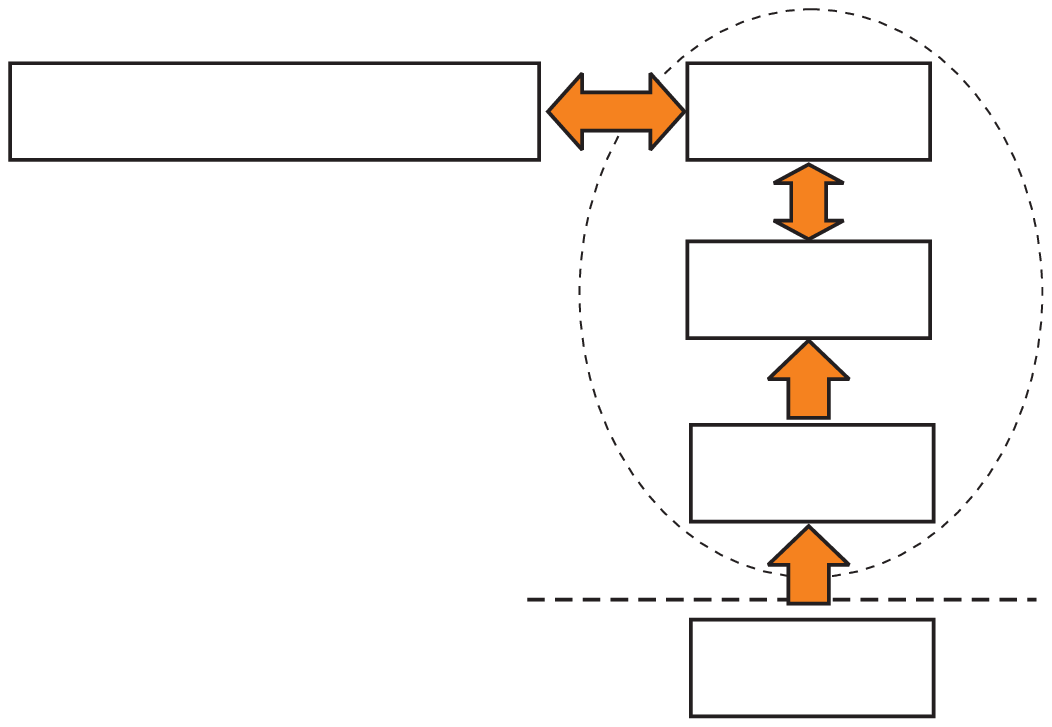}
    \put(20,35){\small{Positivstellensatz}}
    \put(28,30){\small{refutations}}
    \put(10,57){MC problem}
    \put(69,57){$\mathcal{X}_1 \neq \emptyset$}
    \put(69,40){$\mathcal{X}_2 \neq \emptyset$}
    \put(69,22){$\mathcal{X}_3 \neq \emptyset$}
    \put(68,4){$\mathcal{X}_T \neq \emptyset$}
  \end{overpic}
  }
  \caption{Set hierarchy}
  \label{Fig:set_hierarchy}
\end{figure}

To see the relationship of the conditions so far, we show the
relationship among the refutation sets in
Figure~\ref{Fig:set_hierarchy}. The set $\mathcal{X}_1$ and
$\mathcal{X}_2$ are the two largest refutation sets and the
computational complexity of searching these sets grows exponentially
with respect to the dimension of the uncertainty. If we limit the
search for refutations to the set $\mathcal{X}_3$, then computational
complexity is reduced and the results are still tighter than the
existing conditions $\mathcal{X}_T$.

\section{Tightness of quadratic certificates}

In the previous section, we showed that the search for
Positivstellensatz refutations for the matrix cube problem may be
limited to polynomials of low degree with no additional conservatism.
We also provided a hierarchical structure of refutation sets. Although
using quadratic certificates for the refutation is only sufficient, it
is computationally tractable and the conditions are tighter than those
of prior work. In this section, we show that searching quadratic
certificates is also necessary in certain cases. First, we study the
case when the number of uncertainty parameters $m$ is less than or
equal to 2.

\begin{thm}
\label{thm:var_2}
Suppose $m\leq 2$. Then $G(\delta)\geq 0$ for all $\delta$ within the
unit cube if and only if there exists $(S_0,S_1,S_2)\in\mathcal{X}_3$.
\end{thm}

\begin{proof}
  Sufficiency is proved in the previous section and we need only prove
  necessity. Suppose $G(\delta)\geq 0$ for all $\delta$ within the
  unit cube, then
  \begin{align*}
    \bmat{H_0 & H_1 & H_2 & 0 \\
          H_1 & H_0 & 0 & H_2 \\
          H_2 & 0 & H_0 & H_1 \\
          0 & H_2 & H_1 & H_0} \geq 0.
  \end{align*}
  from Lemma~\ref{lem:matrix_cube_LMI}. Without loss of generality,
  assume $H_0$ is nonsingular. Applying the Schur complement gives
  \begin{align*}
    \bmat{H_0 & H_1 & H_2 \\ H_1 & H_0 & 0 \\ H_2 & 0 & H_1} -
    \bmat{ 0 & 0 & 0 \\ 0 & Q_1 & Q_2 \\ 0 & Q_2^T & Q_3} &\geq 0, &
    \bmat{H_0 & H_1 & H_2 \\ H_1 & 0 & 0 \\ H_2 & 0 & 0 } +
    \bmat{0 & 0 & 0 \\0 & Q_3 & Q_2^T \\ 0 & Q_2 & Q_1} &\geq 0
  \end{align*}
  where
  \begin{align*}
    \bmat{Q_1 & Q_2 \\ Q_2^T & Q_3} &= \bmat{H_2 \\ H_1} H_0^{-1}
                                       \bmat{H_2 \\ H_1}^T.
  \end{align*}
  Summing the above two inequalities gives
  \begin{align*}
    \bmat{2H_0& 2H_1 & 2H_2 \\ 2H_1 & H_0 & 0 \\ 2H_2 & 0 &
      H_0} +
    \bmat{0 & 0 & 0 \\
      0 & Q_3-Q_1 & Q_2^T-Q_2 \\ 0 & Q_2-Q_2^T & Q_1-Q_3} \geq 0
  \end{align*}
  Let $W_1=\frac{1}{4}(H_0+Q_3-Q_1)$, $W_2=\frac{1}{4}(H_0-Q_3+Q_1)$.
  Define $J$ as follows
  \begin{align*}
    J = \bmat{H_0-W_1-W_2 & \frac{1}{2}H_1 & \frac{1}{2} H_2 \\[1mm]
              \frac{1}{2} H_1 & W_1 & \frac{1}{4}(Q_2^T-Q_2) \\[1mm]
              \frac{1}{2} H_2 & \frac{1}{4}(Q_2-Q_2^T) & W_2}
  \end{align*}
  and the above inequality implies $J\geq 0$. Let
  \begin{align*}
    S_0 &= \bmat{I \\ \delta\otimes I}^T J \bmat{I \\ \delta\otimes I}, &
    S_1 &= W_1, & S_2 &= W_2.
  \end{align*}
  Expanding $S_0$ shows $(S_0,S_1,S_2)\in\mathcal{X}_3$.
\end{proof}

When $m>2$, if $H_1,\dotsc,H_m$ are either positive or negative
semidefinite, quadratic certificates again provide a necessary
condition.

\begin{thm}
  Suppose the matrices $H_1,\dotsc,H_m$ are either all positive
  semidefinite or all negative semidefinite, then $G(\delta)\geq 0$ for
  all $\delta$ within the unit cube if and only if there exist
  $(S_0,S_1,\dotsc,S_m) \in \mathcal{X}_3$.
\end{thm}

\begin{proof}
  Assume the matrix cube problem holds, then
  \begin{align*}
    x^T H_0 x - \sum_{i=1}^m |x^T H_i x| \geq 0
  \end{align*}
  for all $x\in\R^n$. It also implies
  \begin{align*}
    x^T H_0 x - \frac{1}{2} \sum_{i=1}^m x^T Q_i x -
    \frac{1}{2} \sum_{i=1}^m (x^T Q_i x)^{-1} (x^T H_i x)^2 \geq 0
  \end{align*}
  where
  \begin{align*}
    Q_i &= \begin{cases}
           \displaystyle H_i \qquad \text{ when } H_i\geq 0 \\
           \displaystyle -H_i \hfill \text{ when } H_i\leq 0
           \end{cases}
  \end{align*}
  Apply the Schur complement gives the following polynomial matrix
  inequality
  \begin{align*}
    L(x) &= \bmat{ 2x^T H_0 x - \sum_{i=1}^m x^T Q_i x 
                              & x^T H_1 x & \cdots & x^T H_m x \\
           x^T H_1 x & x^T Q_1 x & & 0 \\
           \vdots & & \ddots & \\
           x^T H_m x & 0 & & x^T Q_m x} \geq 0.
  \end{align*}
  It can be verified that $L(x)$ is SOS since
  \begin{align*}
    L(x) &= e_1 \bbl(2x^T \bbl(H_0-\sum_{i=1}^m Q_i\bbr)x \bbr) e_1^T
       + \sum_{i=1}^m \bmat{e_1 & e_i}
	                 \bmat{x^T Q_i x & x^T H_i x 
			    \\ x^T H_i x & x^T Q_i x}
			 \bmat{e_1 & e_i}^T
  \end{align*}
  Let $\zeta=\bmat{1 & \delta^T}^T$ and
  \begin{align*}
    S_0 &= H_0+\sum_{i=1}^m \delta_i H_i -\frac{1}{2}\sum_{i=1}^m 
                      (1-\delta_i^2) Q_i \\
    S_i &= \frac{1}{2} Q_i \qquad \text{for } i=1,\dotsc,m
  \end{align*}
  Thus, $S_0$ is an sos polynomial matrix since $x^T S_0 x=\zeta^T
  L(x) \zeta$. This shows that $(S_0,S_1,\dotsc,S_m) \in
  \mathcal{X}_3$.
\end{proof}

The above theorems allow us to limit the search to quadratic
certificates in many cases. In general, although the refutation with
quadratic certificates is not necessary and sufficient, it is tight in
many cases. Sometimes, we are also interested in which $x$ and
$\delta$ which minimize $x^T G(\delta) x$. For example, we would like
to know the uncertainty parameter and the state which may almost
destabilize quadratic stability. To investigate this, we examine the
dual of~\eqref{eqn:quadratic_refutation} as follows.
\begin{align*}
    \min \qquad & \trace(H L) \\
    \text{such that } & 
        L = \bmat{L_{00}  & \cdots & L_{0m} \\
                  \vdots & \ddots & \vdots \\
                  L_{m0} & \cdots  & L_{mm}}\geq 0 \\
                {} & L_{ij} = L_{ji} \qquad \text{ for } i,j=0,\dotsc,m \\
                {} & L_{ii} = L_{00} \qquad \text{ for } i=1,\dotsc,m \\
                {} & \trace(L_{00}) = 1
\end{align*}
where
\[
    H = \bmat{H_0 & \frac{1}{2} H_1 & \cdots & \frac{1}{2} H_m \\
               \frac{1}{2} H_1 & 0 & \cdots & 0 \\
               \vdots & & \ddots & \vdots \\
               \frac{1}{2} H_m & 0 & \cdots & 0}.
\]
If there exists a rank one optimal solution $L^{\star}$ and the
optimal value $d^{\star}$ is nonnegative, this implies $G(\delta)\geq
0$ for all $\delta$ within the unit cube. To see this, we can
decompose the rank one solution $L^{\star}$ as
\begin{align*}
  L^{\star} &= \bmat{1 \\ \delta} \bmat{1 \\ \delta}^T \otimes x x^T
\end{align*}
where $\delta\in\{-1,1\}^m$, $x\in\R^n$. These are then the the
optimal solution to the optimization problem
\begin{align*}
  \min_{\substack{\delta\in\{-1,1\}^m \\ x\in\R^n}} \qquad & 
           x^T \bbl( H_0+\sum_{i=1}^m \delta_i H_i \bbr) x.
\end{align*}
Thus this can be interpreted as a lifting of the above dual
optimization problem to a higher dimensional space. This approach
gives a lower bound on the matrix cube problem.

\begin{thm}
\label{thm:sos_exact} 
Suppose $L^{\star}$ is the optimal solution of the dual problem and
$L^{\star}$ is partitioned as above. If $\rank(L^{\star}) =
rank(L_{00}^{\star})$, then $G(\delta)\geq 0$ for all $\delta$ within
the unit cube if and only if the optimal value $d^{\star}$ is
nonnegative.
\end{thm}

\begin{proof}
  We only need to show there is a rank 1 solution of the dual problem.
  Assume $L^{\star}$ is an optimal solution of the dual problem and
  $\rank(L^{\star})= rank(L_{00}^{\star})$. Let $p$ be the rank of
  $\rank(L^{\star})$ and $\rank(L_{00}^{\star})$. We decompose
  $L^{\star}$ as follows,
  \begin{align*}
    L^{\star} &= \bmat{Z \\ & Z \\ & & \ddots \\ & & & Z}
                 \bmat{Y_0 \\ Y_1 \\ \vdots \\ Y_m}
                 \bmat{Y_0 \\ Y_1 \\ \vdots \\ Y_m}^T
                 \bmat{Z \\ & Z \\ & & \ddots \\ & & & Z}^T
  \end{align*}
  where $Z\in\R^{n\times p}$, and $Y_0,\dotsc,Y_m\in\R^{p\times p}$
  have full rank. Since $Z Y_i Y_0^T Z^T = L_{0i}^{\star}\in\S^n$ and
  $Y_0$ has full rank, there exists $D_i\in\S^p$ such that $Y_i=Y_0
  D_i$.  For all $i,j>1$, $i\ne j$, we also know
  \begin{align*}
    L_{ij}^{\star} &= Z Y_i Y_j^T Z^T 
                 = Z Y_0 D_i D_j Y_0^T Z^T \\
                 &= Z Y_0 D_j D_i Y_0^T Z^T 
                 = L_{ji}^{\star}
  \end{align*}
  and the above condition holds if and only if $D_i D_j=D_j D_i$. It
  means that $\{D_1,\dotsc,D_m\}$ is a commuting family and there must
  exist an unitary matrix $U$ which diagonalized $D_1,\dotsc,D_m$.
  Let $D_i=U\Sigma_i U^T$ for $i=1,\dotsc,m$ and $\Sigma_i$ is
  diagonal.  From $L_{ii}^{\star}=L_{00}^{\star}$, we can show
  \begin{align*}
    L_{ii}^{\star} &= Z Y_i Y_i^T Z^T \\
            &= Z Y_0 D_i^2 Y_0^T Z^T \\
            &= Z Y_0 U \Sigma_i^2 U^T Y_0^T Z^T \\
            &= L_{00}^{\star} \\
            &= Z Y_0 Y_0^T Z^T
  \end{align*}
  and the above condition holds if and only if $\Sigma_i^2=1$. Let
  \begin{align*}
    X &= Z Y_0 U = \bmat{x_0 & x_1 & \cdots x_p}, \\
    \Sigma_i &= \bmat{\lambda_{i1} \\ & \ddots \\ & & \lambda_{ip}}
    \quad \text{ for } i=1,\dotsc,m
  \end{align*}
  where $x_i\in\R^{n\times 1}$ and $\lambda_{ij}^2=1$ for
  $i=1,\dotsc,m$, $j=1,\dotsc,p$. We can further decompose $L^{\star}$
  as
  \begin{align*}
    L^{\star} &= \bmat{X \\ X\Sigma_1 \\ \vdots \\ X\Sigma_m}
                 \bmat{X \\ X\Sigma_1 \\ \vdots \\ X\Sigma_m}^T
               = \sum_{k=1}^m 
               \bmat{x_k \\ x_k\lambda_{1k} \\ \vdots \\ x\lambda_{mk}}
               \bmat{x_k \\ x_k\lambda_{1k} \\ \vdots \\ x\lambda_{mk}}^T
              = \sum_{k=1}^m \alpha_k J_k
  \end{align*}
  where $\alpha_k=\norm{x_k}^2$ and
  \begin{align*}
    J_k &= \frac{1}{\norm{x_k}^2} 
          \bmat{x_k \\ x_k\lambda_{1k} \\ \vdots \\ x\lambda_{mk}}
          \bmat{x_k \\ x_k\lambda_{1k} \\ \vdots \\ x\lambda_{mk}}^T
  \end{align*}
  Substituting $L^{\star}$ into the objective function gives
  \begin{align*}
    \trace(L^{\star} H) &= \sum_{k=1}^p \bbl(\trace(J_k H) \bbr)
  \end{align*}
  Note that each $J_k$ is also a feasible solution of the dual problem
  and the optimal value is the linear combination of several feasible
  solutions. Hence, every $J_k$ must also be an optimal solution and
  it is rank 1.
\end{proof}

This theorem gives a condition for existence of a rank one solution.
In addition, if the condition is satisfied, the above construction
also gives a procedure to extract the uncertainty variable $\delta$
which has the smallest eigenvalue in the matrix cube problem.

\section{Examples}

\subsection{Quadratic stability}
\label{sec:exampleone}

\begin{figure}[htb]
  \centerline{\includegraphics[width=0.65\textwidth]{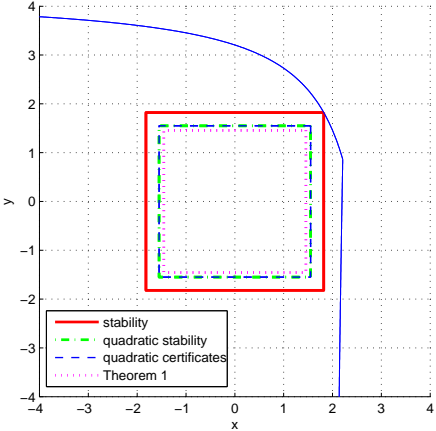}}
  \caption{stability bound from various conditions}
  \label{fig:twod_example}
\end{figure}

In this section, we check quadratic stability of a linear
time-invariant system to provide a specific numerical example. Suppose
the linear time-invariant system with uncertainties is
\begin{align*}
    \dot x &= (A_0 + \sum_{i=1}^m \delta_i A_i) x
\end{align*}
where $\norm{\delta}_{\infty}\leq R$. We would like to compute the
largest $R$ such that the system is quadratically stable for all
$\delta$ within the cube. We first study the two-variable system as
follows.
\begin{align*}
  \dot{x}(t) =
  \bmat{
    -0.4 & 0 & -0.3\delta_1+1 \\
    0 & -3.2 & 0.3\delta_1-0.5 \\
    0.4\delta_2-0.8 & 0.3\delta_2-2.2 & \delta_1-1.7 }
  x(t).
\end{align*}
We compute the largest such $R$ for which all $\delta$ within the
corresponding cube lead to stability.  We similarly compute the
largest cube admitting quadratic stability, and the bounds on this
cube obtained using quadratic certificates and
Theorem~\ref{thm:thm_bental} from Ben-Tal and Nemirovski. These are
shown in Figure~\ref{fig:twod_example}. As discussed in the previous
section, the bound obtained using quadratic certificates is exact, and
we do not need to pursue higher degree certificates.

\subsubsection{general case}

\begin{figure}[htb]
  \centerline{\includegraphics[width=0.55\textwidth]{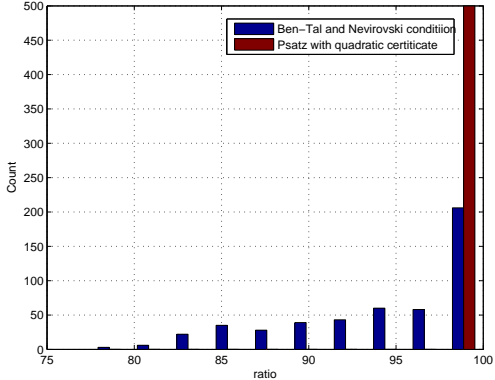}}
  \caption{stability bound from various conditions}
  \label{fig:stability_bound}
\end{figure}

We also test random cases when $m$ is greater than 2.  We let
$x\in\R^5$, $A_0=-4I$ and we choose $m$ ranging from 3 to 8.  For each
$m$, we randomly generate 100 different $A_1,\dotsc,A_m$ and compute
the actual quadratic stability bound $R_e$, the bound $R_s$ from
Positivstellensatz refutation, and the bound $R_t$ using conditions in
Theorem~\ref{thm:thm_bental}. After computing these three values, we
calculate the ratio between $R_s$, $R_e$ and $R_t$, $R_e$ shown in
Figure~\ref{fig:stability_bound}. The bound using quadratic
certificates is very close to the actual quadratic stability bound in
most of the case, but the stability bound from Ben-Tal and
Nemirovski's condition can be $20$ percent smaller than the exact
bound. In these examples, the condition in
Theorem~\ref{thm:thm_positive} indeed gives a tighter stability bound.

\subsection{MAXCUT problem}

\begin{figure}[htb]
  \centerline{\includegraphics[width=0.6\textwidth]{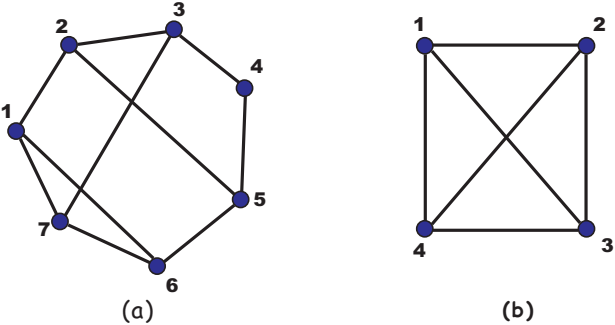}}
  \caption{MAXCUT example}
  \label{fig:maxcut_example}
\end{figure}

Next, we consider the max cut problem. It is well-known that the max
cut problem is an NP hard problem and the SDP relaxation provides a
tractable approach which has 87 percents performance
guarantee~\cite{Goemans_97,Goemans_Williamson_95}. In this example, we
convert the max cut problem to the matrix cube problem and use the
graphs in Figure~\ref{fig:maxcut_example} as examples. 

By using Positivstellensatz refutation with quadratic certificates,
the optimal capacity of the graph in Figure~\ref{fig:maxcut_example}
(a) is $9$ and the rank of the dual problem is 2. We also can infer
the solution is exact from Theorem~\ref{thm:sos_exact} and the optimal
cuts are as follows.
\[
    \begin{tabular}{c|cc}
        node &  \\ \hline
        1 & 1 & -1 \\ \hline
        2 & -1 & 1 \\ \hline
        3 & 1 & -1 \\ \hline
        4 & -1 & 1 \\ \hline
        5 & -1 & 1 \\ \hline
        6 & 1 & -1 \\ \hline
        7 & 1 & -1
    \end{tabular}
\]

We also compute the optimal capacity of the graph in
Figure~\ref{fig:maxcut_example} (b). The optimal value using quadratic
certificate is 4 and the rank of the dual optimal solution is 6.
Although the rank of the dual problem does not satisfy
Theorem~\ref{thm:sos_exact}, it can be verified that the optimal value
is exact. Furthermore, there are 6 different choices which can be
extracted from the dual optimal solution.
\[
    \begin{tabular}{c|cccccc}
        node & & & & & & \\ \hline
        1 & 1 & 1 & 1 & -1 & -1 & -1 \\ \hline
        2 & -1 & 1 & -1 & 1 & -1 & 1 \\ \hline
        3 & -1 & -1 & 1 & 1 & 1 & -1 \\ \hline
        4 & 1 & -1 & -1 & -1 & 1 & 1
    \end{tabular}
\]
To compare the numerical result, we list the capacity bound from the
SDP relaxation and from the Ben-Tal and Nevirovski's condition. The
refutation using quadratic certificates indeed gives a tight result.
\[
    \begin{tabular}{c|cccc}
     & exact capacity & psatz refutation & Ben-Tal & relaxed SDP \\ \hline
        graph 1 & 9 & 9.0000 & 9.925 & 9.8000 \\
        graph 2 & 4 & 4.0000 & 4.0625 & 4.0000
    \end{tabular}
\]

Although the refutation with quadratic certificates usually gives a
tight result, the condition is only sufficient in general. To give a
counterexample, we consider a fully connected graph with 5 nodes. The
maximum capacity of this graph is $6$ and the capacity computed using
quadratic certificates is $6.25$ which is clearly not tight.

\section{Conclusions}

The question of degree bounds for positivstellensatz refutations is
one of significant importance for practical use of semidefinite
programming for matrix polynomial optimization.  In this paper, we
showed that meaningful bounds can be obtained.  We used matrix
Positivstellensatz refutations to test positive semidefiniteness of an
affine function over a given uncertainty set. When the uncertainty set
is a hypercube, we show that the highest degree certificate needed is
$2m$.  Although the certificates are degree bounded, computational
complexity is still high in general. To reduce the complexity, we
study the case of quadratic certificates and show that the bounds
obtained are still tighter than those obtained from existing
conditions. We also show several cases when refutation using quadratic
certificates is exact.  This result may be useful in analyzing and
synthesizing a robust controller for systems with uncertainties and
robust quadratic optimization.

\bibliography{matrix_cube}

\end{document}